\documentclass[12pt]{article}
\usepackage[curve]{xypic}
 \usepackage[usenames,dvipsnames]{pstricks}
\usepackage{epsfig}
\usepackage[active]{srcltx}
\usepackage{amsthm}
\usepackage{amsfonts}
\usepackage{amscd}
\usepackage{amsmath}
\usepackage{amssymb}
\usepackage{color}
\usepackage{xcolor}
\usepackage{array}
\usepackage{epsfig}
\usepackage{graphics}
\usepackage{graphicx}
\usepackage{verbatim}
\usepackage[utf8]{inputenc}

\usepackage{geometry}
    \def\beginp{ {\noindent\em Proof.} }

\def\AAA{{\cal A}}
\def\BBB{{\cal B}}
\def\CCC{{\cal C}}
\def\DDD{{\cal D}}

\def\FFF{{\cal F}}

\def\NNN{{\cal N}}

\def\RRR{{\cal R}}

\def\al{\alpha}

\def\g{\gamma}
\def\G{\Gamma}

\def\ep{\varepsilon}

\def\La{\Lambda}
\def\De{\Delta}
\def\de{\delta}

\def\smm{\smallsetminus}

\def\sgn{\,\text{sgn}}

\def\vecb0{\vec{\bf 0}}
\def\wh{\widehat}

\def\R{\mbox{$\mathbb R$}}

\def\C{\mbox{$\mathbb C$}}

\def\D{\mbox{$\mathbb D$}}
\def\mystrut{{\rule[-2ex]{0ex}{4.5ex}{}}}
\def\ld{ \left|\,\begin{matrix} }
 \def\rd{\end{matrix}\,\right|}

\def\rr{\rightarrow}

 \def\lv{ \left(\begin{matrix} }
 \def\rv{\end{matrix}\right)}

\def\ds{\displaystyle }

\newcommand{\Ref}[1]{ (\ref{#1}) }
\newtheorem{newthm}{Theorem}
\newtheorem{theorem}{Theorem}[section]
\newtheorem{definition}{Definition}[section]
\newtheorem{lemma}[theorem]{Lemma}
\newtheorem{proposition}[theorem]{Proposition}
\newtheorem{corollary}[theorem]{Corollary}

\newcommand{\REFEQN}[1] { \begin{equation}\label{#1} }
\newcommand{\ENDEQN}{\end{equation}}
\newcommand{\REFTHM}[1] { \begin{theorem}\label{#1} }
\newcommand{\ENDTHM}{\end{theorem}}
\newcommand{\REFDEF}[1] { \begin{definition}\label{#1} }
\newcommand{\ENDDEF}{\end{definition}}
\newcommand{\REFNTH}[1] { \begin{newthm}\label{#1} }
\newcommand{\ENDNTH}{\end{newthm}}
\newcommand{\REFPROP}[1]{\begin{proposition}\label{#1} }
\newcommand{\ENDPROP}{\end{proposition} }
\newcommand{\REFLEM}[1]{\begin{lemma}\label{#1} }
\newcommand{\ENDLEM}{\end{lemma} }
\newcommand{\REFCOR}[1]{\begin{corollary}\label{#1} }
\newcommand{\ENDCOR}{\end{corollary} }

\def\bv{{\bf v}}

\def\whx{\widehat{x}}

\def\whf{\widehat{f}}
\def\whI{\widehat{I}}

\def\whc{\widehat{c}}
\def\whu{\widehat{u}}
\def\whv{\widehat{v}}

\def\wtc{\widetilde{c}}

\def\wtI{\widetilde{I}}

\def\wtS{\widetilde{S}}

\def\wt{\widetilde}
\def\wtf{\widetilde{f}}
\def\wtZ{\widetilde{Z}}
\def\Int{\rm Int}

\def\+{\text{\rm\tiny +}}
\def\-{\text{\small -}}

\newcommand{\eng}[1]{ \langle #1  \rangle} 

\begin{document}
\title{Kneading with weights}
\author{  H.H. Rugh and  Tan Lei }
\maketitle

\begin{abstract} We generalise Milnor-Thurston's kneading theory to the setting of piecewise continuous and monotone interval maps with weights.
We define a weighted kneading determinant  $\DDD(t)$ and establish combinatorially two kneading identities, one with the cutting invariant and one
with the dynamical zeta function. For the pressure  $\log \rho_1$ of the weighted system, playing the role of entropy, we prove
that  $\DDD(t)$ is non-zero when $|t|<1/\rho_1$ and has a zero at $1/\rho_1$. Furthermore, our map is semi-conjugate to  an analytic family $h_t, 0<t<1/\rho_1$ of  Cantor PL maps  converging to an interval  PL map $h_{1/\rho_1}$ with equal pressure\footnote{2010 MSC 37E05, 37E15, 37C30 }.
\end{abstract}

\section{Introduction}

 Let $I=[a,b]$. Let $a=c_0<c_1<\cdots < c_{\ell +1}=b$. 
Set $S=\{0,1,\cdots, \ell\}$.   For each $i\in S$, set $I_i=]c_{i}, c_{i+1}[$ and let 
 $f_i: I_i\to  I$ be  a strictly monotone
 continuous map extending continuously to the closure, 
and finally   assign a constant weight $g_i\in \C$.
 
  We say that $(I_i, f_i, g_i)_{i\in S}$ is a {\bf weighted system}. In the particular case that each $g_i$ equals  1, we say also that the
  system is {\bf unweighted}. 
 
 Milnor-Thurston \cite{MT} developed a widely used kneading theory on unweighted systems so that the maps $f_i$ glue together to a  single continuous map $f$.
 Let us recall a list of their results (see also \cite{Ha} for an enlightening introduction to the subject).

Milnor-Thurston introduce a  
 power series matrix $\NNN(t)$, called the {\em kneading matrix},  which records combinatorially the forward orbits of the cutting points. They establish two  identities:

1. The {\bf Main Kneading Identity}, relating  $\NNN(t)$ to the growth of the cutting points of $f^n$ on any subinterval $J$, and taking  the  form 
\REFEQN{main-simple} \g_J (t)\cdot  \NNN(t)= \text{terms involving boundaries of }J;\ENDEQN 

2. The {\bf zeta-function identity},   relating $\NNN(t)$ to a dynamical Artin-Mazur zeta function
 that counts the global growth of the $f^n$-fixed points, taking the  form 
\REFEQN{second}\zeta(t)\cdot \det  \NNN(t)=1.\ENDEQN
 
Using these identities, Milnor-Thurston derive the following important consequences:
 
 3. For $\log s$ the topological entropy of the map, the matrix $\NNN(t)$ is invertible when $|t|<1/s$. If  $s>1$ the matrix $\NNN(t)$ is singular at  $t=1/s$ and the growth rate of the periodic points is precisely $s$.
 
4. If $s>1$, the map is semi-conjugate to a simple model dynamical system which is a continuous PL  (i.e. piecewise-linear) map of slope $s$.

 Most of this theory has been extended by Preston \cite{Pr} to the general unweighted setting without the assumption of global continuity. An advantage to allow discontinuity 
 at the cutting points is that one can treat tree and graph maps as interval unweighted systems after edge concatenation. See for example Tiozzo \cite{Ti}.  There exist also works that treat tree maps as they are.
  See for example Alves-SousaRamos, Baillif and Baillif-deCarvalho  \cite{AS, Ba, BC}. 
  
  An essential 
 difference in Preston's approach as compared to Milnor-Thurston's lies in the proof of the zeta-function identity. Preston's method is purely combinatorial whereas the original proof tests on a concrete example and then studies behaviours under perturbations.
 
In this work we will generalise all four results above to weighted systems, where the pressure
$\log \rho_1$ will play the role of entropy. Points 1-4 will become Theorems \ref{main}, \ref{zeta}, \ref{invertible}, \ref{mainconjugacy} below.

Our setting is identical to that of Baladi-Ruelle \cite{BR}. In their work  they define a weighted  kneading matrix $\BBB$ and a weighted zeta function, and establish a version  the zeta-function identity using a perturbative method similar to that of Milnor-Thurston. For our purpose we will define a somewhat  different  kneading matrix $\RRR$.  

We will not rely on previous established results but instead provide self-contained proofs. In a way
our results recover partially results in \cite{BR,MT,Pr}.

Our proofs will be fairly elementary, with, as the only background, some basic knowledge of complex analysis. The rest is to play carefully with the combinatorics of iterations, following mostly Milnor-Thurston. 

There is however a notable exception, which is about the proof of the zeta-function identity. For this we  choose to follow the combinatorial method of Preston,
along with several significant differences.
Preston  cuts off  the graph above the diagonal in order to count the intersections, 
instead we keep the graph intact but change signs across the diagonal. Preston's kneading matrix is similar to that of Milnor-Thurston, 
by recording the sequence of visited intervals of a critical orbit. Instead we take the kneading matrix $\BBB$ of Baladi-Ruelle, which records the orbit's position relative to every given critical point. We then add one more dimension to $\BBB$ to obtain our kneading matrix $\RRR$, by incorporating the influence of the boundary cutting points (with a somewhat different choice of sign). These modifications are designed to simplify, even in the unweighted case, Preston's proof of the zeta-function identity.
 Preston's idea is to express $-\Big(\log \zeta(t)\Big)'_t$ as the trace of a certain matrix $\FFF$, and then use
repeatedly the Main Kneading Identity to connect $\FFF$ with  the derivative of the kneading matrix. Here,  many choices
are possible  but most give rise to additional correcting terms. Having tested various possibilities we came up with the current choice of the kneading matrix $\RRR$ and a matrix $\FFF$ for which we have the simplest relation  possible,
 i.e.\  $\FFF\RRR=\RRR'$ (see Theorem \ref{fast}). Once this relation established, the zeta-function identity is a one-line computation:
$$-\dfrac d{dt}\log \zeta(t)=Tr \FFF=Tr \RRR'\RRR^{-1}=\dfrac d{dt}\log \det \RRR .$$

The kneading matrix and its smallest positive zero cost relatively little to evaluate. This  enables a fast and accurate computation of the pressure/entropy as well as the semi-conjugacy and the PL model map. 

While experimenting these ideas we noticed that  the system is also semi-conjugate
to a PL map for every $0<t<\rho_1$, although the conjugated system acts on a Cantor set instead of an interval. This numerical observation can  easily be  proved and
has now become our Theorem \ref{conjugacy}. To the best of our knowledge this statement is new, also in the unweighted setting, even though its proof
does not require any new ideas. 

A further justification of our choice of the kneading determinant $\det \RRR$ as compared to  $\det \BBB$, is that the latter  may have a spurious small zero unrelated to the pressure (in Appendix \ref{example} we give an example).

Another originality of this work is the systematic treatment of point-germs relative to points. Each point $x$ in the interior of the interval generates two point-germs: $x^\+$ and $x^\-$. They have often distinct dynamical behaviour and it is convenient to treat the two germs independently. The idea is certainly present to all the papers in the theory.
But highlighting the notion transforms our computations in more concise forms.

Why adding weights to piecewise continuous and monotone maps? One motivation is that one can prescribe slope ratios for the PL model maps, the other is that one can  choose to ignore some parts of a dynamical system by
assigning zero weights, so to reveal deeper entropies hidden for example in renormalisation pieces. 

A further application, not pursued in the current work, is to construct various invariant measures by playing with weights and following Preston's
construction of   measures maximizing the entropy.

{\bf Acknowledgement.} This note originates from the second author's  lecture notes for the ANR LAMBDA meeting in April 2014.
 organised by R. Dujardin. We would also like to thank G. Tiozzo for enlightening discussions.

  \section{Notation and results}

 Let $I=[a,b]$. Let $a=c_0<c_1<\cdots < c_{\ell +1}=b$. 
  Set $I_i=]c_{i}, c_{i+1}[$ and let 
 $f_i: I_i\to  I$ be strictly monotone
 continuous maps for $i=0,\cdots, {\ell }$. We write
$f=(f_0|_{I_0},\cdots, f_{\ell}|_{I_\ell})$ and let
$s_i={\rm sign} ( f_i(c_{i+1}^-) - f_i(c_u^+)) =\pm 1$
denote the sign of monotonicity.
We consider $f$ as undefined at the cutting points.
On the other hand, each $f_i$ extends to a continuous map from the closed interval
$[c_{i}, c_{i+1}]$ to $[a,b]$.

We call $\CCC(f)=\{c_i : 1\leq i\leq \ell\}$ the interior
 cutting points of the interval. 
 The set of cutting points, $\CCC_*(f)$, includes $c_0$ and $c_{\ell+1}$.

In order to treat monotonicity and discontinuities in a consistent manner
it is convenient to extend our base interval $I$ to its unit-tangent
bundle, also denoted the space of point-germs
$\wh I $: 
each point $x\in I\smm\{a,b\}$ generates two point-germs denoted
$x^\text{\rm\tiny +}=(x,+1)$ and $x^\text{\small -}=(x,-1)$ while
the boundary points $a, b$ each has only one point-germ
$a^\text{\rm\tiny +}$ and $b^\text{\small -}$.
We write   $$  \ep(x^\text{\rm\tiny +}):=1\text{\ \  and\ \ \ }
 \ep(x^\text{\small -}):=-1$$ for the direction of the germ.
In order to make some formulae in Section \ref{zeta-identity} more
concise, we  set (artificially)  {\bf $c_0^\-=b^\-$}
so that $\{c_0^\+, c_0^\-\}=\{a^\+, b^\-\}$. 
For $x\in I$ we denote by
$\whx=(x,\sigma)$
the point-germ based at $x$ and in the direction 
$\sigma\in \{\pm 1\}$.

It is notationally 
convenient to define an order $<$ on the collection of point-germs 
together with base points,  $I\cup\wh I$,
by declaring that for two base points $x<y$ we have
$x<x^\text{\rm\tiny +}<y^\text{\small -}<y< y^\+$. 
Given two point-germs $\wh u, \wh v\in \wh I$ with $\wh u<\wh v$,
we define  $$\eng{\wh u, \wh v}:= \Big\{x\in I\mid  \wh u< x< \wh v\Big\}$$
as a sub interval of $I$. 
It is then consistent
to write e.g.\  $[u,v[=\eng{u^\text{\small -}, v^\text{\small -}}$ and
 $]u,v[=\eng{u^\text{\tiny +},v^\text{\small -}}$. 
 Note that the boundary points $a,b$ never belong to an interval of 
 the form $\eng{\wh u, \wh v}$. When $J=]u,v[$ is an open interval
 we set $\wh J = \{ \whx : u<x<v \} \cup \{u^+\} \cup \{v^-\}$.
 In particular, 
 $\wh I_i=\{ \whx : c_i<x<c_{i+1} \} \cup \{c_i^+\} \cup \{c_{i+1}^-\}$,
 $0\leq i \leq \ell$. We observe that the $\wh I_i$'s are disjoint and
 their union is $\wh I$.

Our original map 
$f$ induces a well-defined map $\whf : \wh I \rr \whI$. 
When $\wh x=(x,\sigma) \in \wh I_i$ then $f(\whx)=(y,\sigma')$ is 
simply the germ based at $y=\lim_{t\rr 0^+} f_i(x+\sigma t)$ whose direction is
$\sigma'= s_i \sigma$. Note that on each $\whI_i$, $\whf$ is monotone
because $f$ is {\em strictly} monotone. We will usually
write $f$ also for the extended map $\whf$.

  For each $0\leq i \leq \ell$ we  let $g_i\in \C$ be a weight
associated with the interval $I_i$. Both $g_i$ and $s_i$ gives rise to
functions on $\wh I_i$ by declaring
$s(\wh x)=s_i$ and $g(\wh x)=g_i$ whenever $\wh x\in \wh I_i$.
We may define products along orbits,
$s^n, g^n, [sg]^n$ by setting $s^0=g^0=1$, and 
$$\forall\,n\ge1,\quad s^n(\whx):=\prod_{k=0}^{n-1} s(f^k(\whx)), \quad g^n(\whx):=\prod_{k=0}^{n-1} g(f^k(\whx)), \quad [sg]^n:=s^ng^n.$$
Note that $s^n(\whx)$ is the sense of monotonicity of $f^n$ at $\whx$.

We define  a half-sign function: $$ \forall\ \whx\in \wh I,\ y\in I,\quad \sigma(\whx,y):=\dfrac12 \sgn(\whx-y)=\left\{\begin{array}{ll} +1/2 & \text{if  $\whx>y$}\\ -1/2 & \text{if $\whx<y$ }
\end{array}\right. .$$

Concerning forward orbits of point-germs we
set for  $j,k=0,\cdots, \ell $ :
 \begin{eqnarray} \theta(\whx,t;c_k)&=& \sum_{m\ge 0}t^m [sg]^m(\whx)\cdot \sigma(f^m\whx, c_k), \text{ so in particular}\\
 \theta(\whx,t;c_0)&=& \sum_{m\geq 0}t^m [sg]^m(\whx)\cdot \sigma(f^m\whx, c_0)=  \dfrac12 \sum_{m\ge 0}t^m [sg]^m(\whx)\\
\ep^*(\wh c) &=& \ep(\wh c) \quad \text{ if } \wh c\ne c_0^\pm\ \ \text{and}\ \ 
\ep^*(\wh c) \ =\  +1\quad  \text{ if } \wh c= c_0^\pm\\
  R_{jk}(t)&=& \sum_{\wh c_j =c_j^\+, c_j^-} \ep^*(\wh c_j)\cdot \theta(\wh c_j,t;c_k),\quad  \\
\RRR(t) &=& (R_{jk}(t) )_{0\leq j,k\leq \ell} \ \text{ (the {\bf kneading matrix})} \\
\BBB(t) &=& (R_{jk}(t) )_{1\leq j,k\leq \ell} \ \text{ (the {\bf reduced kneading matrix})} 
\end{eqnarray}\label{theta}
In particular, one has (note the signs):
\begin{eqnarray}
R_{{ j}k}(t) &=& \theta( c_j^\text{\rm\tiny +},t;c_k){ -}\theta(c_j^\text{\small -},t;c_k) =:\De_{c_j}\theta(\cdot,t;c_k) \ { (j>0)}  \ \text{while}\\
  R_{{ 0}k}(t)&=& \theta(a^\+,t;c_k){ +}\theta(b^\-,t;c_k)
\end{eqnarray}
(this  choice of signs is designed to absorb boundary correcting terms in later calculations).

Regarding 'backward'-orbits we
define $Z_1$ as the set of  level-$1$ cylinders 
$(j):=I_j=]c_j, c_{j+1}[$, $j=0,1,\cdots,\ell$.
Define then recursively $Z_n$ as the set of non-empty
level-$n$  cylinders of the form~:
$(i_0i_1\cdots i_{n-1}):=I_{i_0}\cap f_{i_0}^{-1}(i_1\cdots i_{n-1})$.
Each $\al=(i_0i_1\cdots i_{n-1})$ is an open interval 
$]u,v[=\eng{u^\text{\tiny +},v^\text{\small -}}$. We
set $\wh \partial \al=\{u^\text{\tiny +},v^\text{\small -}\}$.
For $0\le j< n$, $f^j(\al)\subset I_{i_j}$. 
So  $f^n$ maps $\al$ homeomorphically onto its image, 
in particular  each of the functions $s^j$ and $g^j$, $0\le j<n$,
is constant on $\al$. 

\REFDEF{expansive}
We call \; $(I_i,f_i)_{0\leq i\leq \ell}$\;  expansive if \
$\ds \lim_{n\rr \infty} \sup_{\al\in Z_n} {\rm diam}\; (\al) = 0$.
\ENDDEF
  
 For any $y\in I$, set   $\G_{0,y}=\{y\}$, and for $p>0$, 
 $$ \G_{p,y}=\Big\{x \in \bigcup_{\al\in Z_p}\al\ \big| \ f^p(x)=y\Big\}.$$
 Note that $x\in \G_{p,y}$ implies that
 $g^p(x^\text{\small -})=g^p(x^\text{\rm\tiny +})$, 
 for which we simply write $g^p(x)$. This is because $g^0(x)\equiv 1$ and
 every $j$-iterate ($0\le j<p$) of a $p$-cylinder $\al\in Z_p$ belongs to  some level-1 cylinder.
   Define $$\g_{y}(t)=\sum_{p\ge 0} \sum_{x\in \G_{p,y}}t^p g^p(x) \quad \text{and, for }J\subset ]a,b[,\quad  \g_{y, J}(t)=\sum_{p\ge 0} \sum_{x\in \G_{p,y}}t^p g^p(x) \chi_J(x).$$

This function counts the (weighted) number\footnote{
In the case $g_i\equiv 1$, we have
$\ds  \g_{y, J}(t)=\sum_{p\ge 0}\#(\G_{p,y}\cap J) \,t^p$.} of preimages of $y$.
   
   Clearly when $J$ and $J'$ are disjoint subsets we have
   $$ \g_{y, J}(t)+ \g_{y, J'}(t)= \g_{y, J\cup J'}(t).$$
\REFTHM{main}(Main Kneading Identity, or MKI in short) 
   For  any interval $J=\eng{\whu,\whv}$ in $I$,
\REFEQN{identity} 
    \forall\ k\in\{0,\cdots, \ell\},\quad  
    \sum_{j=1 } ^\ell \g_{c_j, J}(t) R_{jk}(t)= 
    \theta(\wh v,t;c_k)-\theta(\wh u,t;c_k) =: \De_J\theta(\cdot,t;c_k)\ENDEQN 
(the term $j=0$ is not included in the sum, but  we do allow $k=0$).
  \ENDTHM

   We also need a particular way to count the fixed points of $f^n$.  
  Fix $n\ge 1$ and an $n$-cylinder $\al$. The value of $g^n(x)$ is a constant on $\al$, denoted by $g^n_{|\al}$. We define the (fixed point counting) weight $\omega(\al)$ by
  $$\omega(\al)=-g^n_{|\al}  \sum_{\whx\in \wh \partial \al } \sigma(f^{n}\whx,x)\cdot \ep(f^n\whx).$$
  We refer to Appendix  \ref{count} for an account of
  the geometric meaning of this weight. We then introduce the following 
  weighted counting of  fixed points of $f^n$:
  $$N_n:=\sum_{\al\in Z_n} \omega(\al).$$
  Set $Z(t):=\ds \exp\Big(\sum_{n\ge 1} \dfrac1n N_nt^n\Big),\quad N_f(t):=\ds \sum_{n\ge 1} N_nt^{n-1}=(\log Z)'.$

\REFTHM{zeta} (zeta-function identity) 
   Let $\DDD(t)=\det \RRR(t)$. We then have $Z(t)\cdot \DDD(t)=1$,
   or equivalently $$N_f(t)+\dfrac{\DDD'(t)}{\DDD(t)}=0.$$\ENDTHM


\REFDEF{pressure}
For every $n\ge 0$ we write $\ds \big\|g^n\big\|_\infty=\sup_{\al\in Z_n} \left|g^n_{| \al}\right|$ and $\ds\big\|g^n\big\|_1=\sum_{\al\in Z_n} \left|g^n_{| \al}\right|$. We then set 
\begin{equation}
\rho_\infty: =\limsup_{n\to \infty}
\big\|g^n\big\|_\infty^{1/n} \quad
\leq \quad
\rho_1: =\limsup_{n\to \infty}
\big\|g^n\big\|_1^{1/n} .
\end{equation}
We also call $\log \rho_1$ the {\bf pressure}\footnote{In the case $g_i\equiv 1$, we have $\rho_\infty=1$ and $\rho_1$  is the growth rate of the $n$-cylinders. By Misiurewicz-Szlenk (\cite{MS})
$\log \rho_1$ is equal to the topological entropy of the unweighted system. } of the
weighted system $(I_i,f_i,g_i)_{i\in S}$. This is consistent
with  usual "thermodynamic formalism" for dynamical systems.
\ENDDEF
 
\REFTHM{invertible}  We have
\begin{enumerate}
\item The power series for
$\theta(\wh x,t;c_k)$, $R_{jk}(t)$ define analytic functions
of $t$ on the disc $\{ |t|<1/\rho_\infty\}$.
\item The kneading matrix $\RRR(t)$ is invertible when $|t|< 1/\rho_1$.
\item Suppose that $\rho_1>\rho_\infty$ and all $g_i\geq 0$. Then
  $\RRR(t)$ is non-invertible at $t=1/\rho_1$ and $1/\rho_1$ coincides with the radius of convergence of $Z(t)$.\\
\end{enumerate}
\ENDTHM

\REFTHM{conjugacy}
 Assume $\rho_1>\rho_\infty$ and all $g_i >   0$. 
  For each $0<t<1/\rho_1$  there is a monotone (non-continuous)
  map $\phi_t : \wh I \rr  [0,1]$ with the following properties
\begin{enumerate}
\item [{\bf A.}] For $0\leq i\leq \ell$, let
 $\wt I_{t,i} = [ \phi_t(c_i^+), \phi_t(c_{i+1}^-)]$
(this is an interval or a point).
 The collection $\wt I_{t,i}$, $0\leq i \leq \ell$
     is pairwise disjoint.
\item [{\bf B.}] 
  For each $i$ there is an affine map
  $\wt f_{t,i} : \wt I_{t,i} \rr [0,1]$ of slope $s_i/(t g_i)$ such that 
   \begin{equation}
         \phi_t \circ f_{\ | \;\whI_i} = 
	   \wt f_{t,i} \circ \phi_t \mbox{}_{\ | \;\whI_i} 
   \end{equation}
\item [{\bf C.}] 
The partially defined dynamical system 
$(\wt I_{t,i},\wt f_{t,i})_{0\leq i\leq \ell}$
is uniformly expanding and its maximal invariant domain is precisely
$\Omega_{t}= \phi_t(\wh I)$.
\end{enumerate}
  \ENDTHM

 \noindent {\bf Remark.} 
  Thus, $\phi_t$ is
  semi-conjugating  the dynamical system $(\wh I, f)$
   to a uniformly expanding dynamical system
  $(\Omega_{t,i}, \wt f_{t,i})_{0\leq i\leq \ell}$.
   In general, $\Omega_t$ is a Cantor set. The subset
   $\Omega_{t,i}=\phi_t(\wh I_i)$ is trivial, i.e. 
   reduced to a point, precisely when the forward orbit of $I_i$ never
   encounters a cutting point. This can not happen if the 
   original system is expansive.
  
 The proof will show that the semi-conjugacy $\phi_t$ can be explicitly expressed as $$\dfrac{h(\whx)-h(a^\+)}{h(b^\-)-h(a^\+)}\quad \text{with}\quad
   h(\whx)=\Big( \theta(\whx,t;c_k)\Big)_{k=0,\cdots \ell} \cdot \RRR^{-1} \cdot \lv 0\\ G(c_1,t) \\ \vdots \\ G(c_\ell, t) \rv,$$
where $G(x,t)$ is the average of the generating functions for $g^n(x^\-)$ and $g^n(x^\+)$ (see \Ref{generating}, \Ref{w-identity} and \Ref{phi-def}).  If $\ell=1$, 
one can replace $h(\whx)$ by $ \theta(\whx,t;c_1)$, which is  particularly simple to implement numerically.

   When taking the limit as
   $t\nearrow 1/\rho_1$ we obtain a different type of semi-conjugacy:

\REFTHM{mainconjugacy} 
  Assume $\rho_1>\rho_\infty$ and all $g_i\geq 0$. 
   There is a monotone continuous surjective map $\phi : I \rr [0,1]$
   with the following properties:
   Denote by  $\wtS\subset S:= \{0,\cdots,\ell\}$ the subset
   of $i$'s for which  
   $\wt I_i= {\rm Int}\; \phi(I_i)$ is non-empty. Then

 \begin{itemize}
 \item[{\bf A.}]
     For every $i\in \wtS$, there is an affine map
     $\wtf_i$ of 
      slope $s_i \rho_1 /  g_i$ such that
      \[ \wtf_i(\phi(x))=\phi(f_i(x)), x\in I_i.\]
 \item[{\bf B.}]
The two weighted systems  $(I_i,  f_i,g_i)_{i\in S}$ 
and $(\wtI_i,  \wtf_i,g_i)_{i\in \wtS}$ have equal pressures.\\
    \item [{\bf C.}]
     If the sytem $f=(I_i,f_i)_{i\in S}$ extends continuously to a map on
     $[a,b]$ then so does $\wtf = (\wtI_i,\wtf_i)_{i\in \wtS}$ 
     on $[0,1]$ and $\phi$ gives a genuine topological semi-conjugacy.
     We have in this case for every $x\in [a,b]$:
      \[ \wtf(\phi(x))=\phi(f(x)).\]
\end{itemize}
\ENDTHM
 
 For the last theorem, 
 some intervals may disappear under the semi-conjugacy, i.e.\ the
 set $\wtS$ becomes a strict subset of $S$.
  This happens
 in particular, when the original system is not transitive and
 contains sub-systems of a smaller pressure. The set  $\wtS$ may even
 depend on the choice of the weights $g_i$. 
 In particular, intervals for which $g_i=0$ disappear under the conjugacy.

   \section{The Main Kneading Identity}

\REFLEM{constant term}
We have $\RRR(0)=id$.\ENDLEM
\beginp
Note that $$R_{jk}(t)=\sum_{\wh c_j=c_j^\text{\rm\tiny +}, c_j^\text{\rm\small -}}\ep^*(\wh c_j)\cdot \theta(\wh c_j,t;c_k)=\sum_{n\ge 0}t^n \sum_{\wh c_j=c_j^\text{\rm\tiny +}, c_j^\text{\rm\small -}}\ep^*(\wh c_j)[sg]^n(\wh c_j) \cdot \sigma(f^n\wh c_j,c_k)$$

 By convention $f^0=id$. Recall that $\ep^*(\wh c_j)=\ep(\wh c_j)$ if $j\ne 0$ and $\ep^*(\wh c_0)=1$.

Assume first $j>0$. Then, for all $k=0,\cdots,\ell$,

$\ds R_{jk}(0)=\sum_{\wh c_j=c_j^\text{\rm\tiny +}, c_j^\text{\rm\small -}}\ep^*(\wh c_j)\cdot [sg]^0(\wh c_j) \cdot \sigma(f^0\wh c_j,c_k)=\sum_{\wh c_j=c_j^\text{\rm\tiny +}, c_j^\text{\rm\small -}}  \ep(\wh c_j)\cdot \sigma( \wh c_j,c_k) =\de_{jk}.$

Also, $\ds R_{0k}(0)=\theta(a^\+,0;c_k)+\theta(b^\-,0;c_k)=
\sigma(a^\+,c_k) + \sigma(b^\-,c_k)= \de_{0k}$.\qed

\subsection{Proof of Theorem \ref{main}.}

Consider first an open interval $J=]u,v[\subset ]a,b[$, and a $c_k$ for some $k\in \{0,\cdots, \ell\}$.
For each $n\ge 0$ and each $(n+1)$-cylinder $\al\in Z_{n+1}$, the functions
$\ds [sg]^n(\whx)=\prod_{j=0}^{n-1} s(f^j\whx)g(f^j\whx)$ and $\sigma(f^n\whx,c_k)$, $\whx\in \wh \al$ are 
constants.
When $\al\in Z_{n+1}$ and $\al\cap J\ne \emptyset$, then obviously
$$\sum_{\whx\in \wh \partial (J\cap \al)}\ep(\whx)= 1+(-1)=0.$$
So the following power series vanishes identically:
$$\sum_{n\ge 0} t^n\sum_{\al\in Z_{n+1},\whx\in \wh \partial (J\cap \al)}\ep(\whx)\cdot [sg]^n(\whx)\cdot \sigma(f^n\whx,c_k)=0.$$
In this sum, $\whx =u^\+, v^\-$ appears for every $n\ge 0$.
Extracting their contributions we write:
\REFEQN{num}
  \sum_{\whx\in \wh \partial J}\theta(\whx,t;c_k)\cdot \ep(\whx)
   + \sum_{n\ge 0} t^n\sum_{\al\in Z_{n+1},\whx\in \wh \partial \al}
 \chi_J(x)\cdot \ep(\whx)\cdot [sg]^n(\whx)\cdot \sigma(f^n\whx,c_k)=0.
\ENDEQN
Now
when $\al\in Z_{n+1},\whx\in \wh \partial (J\cap \al)$, there is a unique minimal integer $0\le p\le n$ for which $f^p(\wh x)=\whc$ for some
$c\in \{c_1,\cdots, c_\ell\}=:\CCC(f)$ and $\whc =c^\+$ or $c^\-$ 
(note that the boundary points $a,b$ are excluded here, 
since for an interior point to be mapped to them,
it has to pass an interior cutting point just before).
Recall that $\G_{p,c}=\{x\in \bigcup _{\al\in Z_p}\al\mid f^px=c\}$
and $\G_{0,c}=\{c\}$. 
When $x\in \G_{p,c}$ and $f^p\whx=\whc$, then
$g^p(\whx)=g^p(x)$, $\sigma(f^p\whx,c_k)=\sigma(f^{n-p}\whc, c_k)$
and also (the essential point here
is that  the sign $s^p(\wh x)$ is absorbed in $\ep(\whc)$)
$$\ep(\whx)\cdot [sg]^n(\whx)=
   g^p(x)\Big(\ep(\whx)s^p(\whx)\Big)[sg]^{n-p}(\whc)=
   g^p(x)\cdot \ep(\whc)\cdot [sg]^{n-p}(\whc).$$
So we obtain, for the second term in \Ref{num} (writing $t^n=t^pt^q$),
$$ \sum_{c\in \CCC(f)}\Big[\Big(\sum_{p\ge 0} t^p\sum_{x\in \G_{p,c}}
     g^p(x) \chi_J(x)\Big)\cdot  
\sum_{\whc=c^\pm,q\ge 0}
     t^q \cdot \ep(\whc)\cdot [sg]^q(\whc)\sigma(f^q\whc,c_k)\Big]$$
$$= \sum_{c\in \CCC(f)}\g_{J,c}(t)\cdot \De_c\theta(\cdot,t;c_k).$$
Combining with \Ref{num} we get the Main kneading Identity when $J$ is an open interval. 

It remains to prove the case that $J$ is half closed or closed. 
Consider for example $J=\eng{u^\text{\rm\small -},v^\text{\rm\small -}}$ with
$a<u<v\leq b$.
We have $\eng{a^\text{\rm\tiny +},v^\text{\rm\small -}}= \eng{a^\text{\rm\tiny +},u^\text{\rm\small -}}\sqcup J$ and the
additivity $\g_{c, \eng{a^\text{\rm\tiny +},v^\text{\rm\small -}}}= \g_{c, \eng{a^\text{\rm\tiny +},u^\text{\rm\small -}}}+ \g_{c, J}$.
The result then follows by applying the identity to the two intervals $\eng{a^\text{\rm\tiny +},u^\text{\rm\small -}}$ and $\eng{a^\text{\rm\tiny +},v^\text{\rm\small -}}$ and
subtracting. 
\qed

\section{Zeta functions and kneading determinants}\label{zeta-identity}

In this section we prove Theorems \ref{zeta} and \ref{invertible}.

\subsection{Relating $N_f(t)$ to $\RRR'(t)$}

Set $\wh \CCC(f):=\{a^\+=c_0^+, c_1^\-, c_1^\+,\cdots, c_\ell^\-, c_\ell^\+, b^\-=c_0^-\}$. For all $\wh c\in \wh\CCC(f)$, 
set $\G_{0,\wh c}=\{\wh c\}$ and, for $p\ge 1$,  $$\G_{p, \wh c}=\{\whx\in \wh I\mid f^p\whx=\whc,\ , f^j \whx \notin \wh\CCC(f) \text{ for }0\le j <p\}.$$

If $\wh c \ne c_0^\pm$, then for any $\whx\in \G_{p,\wh c}$ we have $x\in \G_{p,c}$. Conversely for any $x\in \G_{p,c}$ 
exactly one of $x^\pm$ belongs to $\G_{p,\wh c}$.

Notice that if $\wh c = c_0^\pm$ then $\G_{p,\wh c}=\emptyset $ when $p\ge 1$: due to the forward invariance of $I$ we have  $f^{-1}(\{a,b\})\subset\{a,b,c_1,\cdots, c_\ell\}$, so every orbit passing though $\{a,b\}$ must pass through $\{c_1,\cdots, c_\ell\}$ just before.

Fix $n\ge 1$ and 
an $n$-cylinder $\al$. 
Note that for each $\whx\in \wh \partial \al$, we have $g^n_{|\al}\cdot  \ep(f^n\whx)=[sg]^n(\whx)\cdot \ep(\whx)$, so
$$- \omega(\al):=g^n_{|\al}  \sum_{\whx\in \wh \partial \al } \sigma(f^{n}\whx,x)\cdot \ep(f^n\whx) =\sum_{\whx\in \wh \partial \al } \sigma(f^{n}\whx,x) [sg]^n(\whx)\cdot \ep(\whx).$$
   To each  $\whx\in \wh \partial \al$, there is a unique $\wh c\in \wh \CCC(f) $
   and $0\le p<n$ such that $x\in \G_{p.\whc }$. Schematically,  
   $$\whx \xrightarrow[p\text{ minimal}]{f^p} \whc \xrightarrow[]{f^{q+1}} f^n\whx=f^{q+1}\whc$$
   
Setting $q$ such that $p+q=n-1$, we have the ``co-cycle'' properties:
$$s^n(\whx)\ep(\whx)=s^{q+1}(\whc)\ep(\wh c),\  g^0(\whx)=1,\ g^n(\whx)=g^{q+1}(\whc )g^p(\whx).$$

Now \begin{eqnarray*}- \sum_{n\ge 1} t^{n-1} N_n&=&
  - \sum_{n\ge 1} t^{n-1}\sum_{\al\in Z_n} \omega(\al)\\ 
  &=& \sum_{n\ge 1} t^{n-1}\sum_{\al\in Z_n,\whx\in \wh \partial \al } 
  \sigma(f^{n}\whx,x) [sg]^n(\whx)\cdot \ep(\whx)\\
  &=&\ds\sum_{\wh c \in \wh\CCC(f)} \sum_{q\ge 0}
  t^q [sg]^{q+1}(\wh c)\cdot \ep(\whc)  
  \sum_{p\ge 0,\,\whx\in \G_{p, \whc}} t^p g^p(\whx)\sigma(f^{q+1} \wh c,x)  \\
&=& \ds\sum_{\wh c \in \wh\CCC(f)\smm\{c_0^\pm\}}
 \sum_{q\ge 0} t^q [sg]^{q+1}(\wh c)\cdot \ep(\whc)
    \Big(\sum_{p\ge 0,\,\whx\in \G_{p, \whc}}
      t^p g^p(\whx)\sigma(f^{q+1} \wh c,x) \Big) \\
 &&+
   \ds\sum_{\wh c=c_0^\pm} \sum_{q\ge 0} t^q [sg]^{q+1}(\wh c)\cdot
   \Big(\ep(\whc) \sum_{p\ge 0,\,\whx\in \G_{p, \whc}}
   t^p g^p(\whx)\sigma(f^{q+1} \wh c,x)  \Big) \end{eqnarray*}
Note that the $\ep(\whc)$ factor in the last expression 
is treated differently for $\whc=c_1^\pm, \cdots, c_\ell^\pm$
and $\wh c=c_0^\pm$. The reason for this is that we want the two expressions
in the parenthesis to be independent of the direction of $\wh c$.
Indeed, for any $\wh u\in \wh I$,
$$\text{for}\quad \whc=c_1^\pm, \cdots, c_\ell^\pm,\quad m_{\wh c}(\whu, t):=\sum_{p\ge 0,\,\whx\in \G_{p, \whc}} t^p g^p(\whx)\sigma(\whu,x)=\sum_{p\ge 0,\,x\in \G_{p, c}} t^p g^p(x)\sigma(\whu,x)$$
$$\text{and for}\quad \whc= c_0^\pm ,\quad   m_{\wh c}(\whu, t):=\ep(\wh c)\sum_{p\ge 0,\,\whx\in \G_{p, \whc}} t^p g^p(\whx)\sigma(\whu,x)=\ep (\wh c)\sigma(\wh u,c) \equiv\dfrac12$$
where we have used the facts that    $\G_{p,c_0^\pm}=\emptyset$ for $p>0$, $g^0(\whx)\equiv 1$ and $g^p(x^\+)=g^p(x^\-)=:g^p(x)$ for $x\in \G_{p,c_j}$, $j>0$.

In both cases $m_{\wh c}(\whu, t)$
is independent of $\ep(\wh c)=+$ or $-$, so we may safely write $m_{c}(\whu, t)$ for this quantity. 
To compactify the two cases we set $\CCC^*(f)=\{c_0, c_1,\cdots, c_\ell\}$.   Recall that $c_0^\+=a^\+$, $c_0^\-=b^\-$ and $ \ep^*(\wh c):=\ep(\wh c)$
if $\wh c\ne c_0^\pm$ and $ \ep^*(\wh c)=1$ otherwise. Then
 \begin{eqnarray*} -\sum_{n\ge 1} t^{n-1} N_n&=&\sum_{c\in \CCC^*(f)}\sum_{\wh c=c^\pm,\,q\ge 0} t^q [sg]^{q+1}(\wh c)\cdot  \ep^*(\whc) \cdot m_c(f^{q+1}\whc,t)\end{eqnarray*}

A central idea (due to Preston) is to consider the right hand side as the 
trace of an $(\ell+1)\times (\ell+1)$ matrix $\FFF$, and to define $\FFF$
in a way so that $\FFF\RRR$ becomes related to $\RRR'$.
  There are many choices suitable for this purpose with most choices
  giving rise to additional correcting terms. There is, however,
  a choice for which the relationship becomes particularly simple
  (note the $*$ in the epsilon factor):

For $i,j\in \{0,1,\cdots,\ell\}$, define
$$\quad F_{ij}(t)=\ds\sum_{q\ge 0,\,\wh c_i=c_i^\pm }
  t^q [sg]^{q+1}(\wh c_i)\cdot \ep^*(\wh c_i)\cdot 
   m_{c_j}(f^{q+1}\wh c_i,t).$$
We then have:

\REFTHM{fast} 
 $\FFF\RRR=\RRR'$.\ENDTHM
 
 \beginp 
We establish at first a consequence of the Main Kneading Identity:

\noindent {\bf Claim.} For  every $\wh w\in \wh I$, $k=0,1,\cdots,\ell$, $$\sum_{j=0}^\ell m_{c_j}(\wh w,t) R_{jk}(t)=\theta(\wh w,t;c_k).$$
\beginp By the Main Kneading Identity, we sum first over interior cutting points:
\begin{eqnarray*}\sum_{j=1}^\ell m_{c_j}(\wh w,t) R_{jk}(t)&=&\sum_{j=1}^\ell \sum_{p\ge 0,\,x\in \G_{p, c_j}} t^p g^p(x)\sigma(\wh w,x)\cdot R_{jk}(t)\\
&=&\sum_{j=1}^\ell \sum_{p\ge 0,\,x\in \G_{p, c_j}} t^p g^p(x)\dfrac12\Big(\chi_{(a^\+,\wh w)}(x)-\chi_{(\wh w, b^\-)}(x) \Big)R_{jk}(t)\\
&=& \dfrac12\Big(2\theta(\wh w,t;c_k)-\theta(a^\+,t;c_k)-\theta(b^\-,t;c_k)\Big)
\end{eqnarray*}
Adding the boundary term
$\ds m_{c_0}(\wh w,t) R_{0k}(t)=
     \dfrac12 \Big(\theta(a^\+,t;c_k)+\theta(b^\-,t;c_k)\Big)$
we get the desired result and end the proof of the claim.

Now, for $i,k\in \{0,\cdots,\ell\}$,
\begin{eqnarray*}   \sum_{j=0}^\ell F_{ij}R_{jk}
 &=& \ds\sum_{q\ge 0,\,\wh c_i=c_i^\pm }
       t^q [sg]^{q+1}(\wh c_i)\cdot \ep^*(\wh c_i)\cdot
   \sum_{j=0}^\ell \Big( 
    m_{c_j}(f^{q+1}\wh c_i,t) R_{jk}\Big)\\
 &=& \ds\sum_{q\ge 0,\,\wh c_i=c_i^\pm} 
       t^q [sg]^{q+1}(\wh c_i)\cdot \ep^*(\wh c_i)  \cdot
         \theta(f^{q+1}(\wh c_i),t;c_k)\\ 
 &=& \sum_{q\ge 0,\,\wh c_i=c_i^\pm}
       t^q [sg]^{q+1}(\wh c_i)\cdot \ep^*(\wh c_i)
       \sum_{p\ge 0} t^p[sg]^p(f^{q+1}\whc_i)\sigma(f^p(f^{q+1}\whc_i),c_k)\\
 &=& \sum_{\wh c_i=c_i^\pm} \sum_{p,q\ge 0}t^{p+q} [sg]^{p+q+1}(\wh c_i)\cdot
     \sigma(f^{p+q+1}\whc_i,c_k)\cdot \ep^*(\wh c_i)\\ 
 &=&\sum_{\wh c_i=c_i^\pm} \Big( \sum_{n\ge 1}n\cdot
     t^{n-1}  [sg]^{n}(\wh c_i)\cdot 
          \sigma(f^{n}\whc_i,c_k)\Big) \ep^*(\wh c_i)\\
 &=& \sum_{\wh c_i=c_i^\pm}
    \Big( \dfrac d{dt}\theta(\wh c_j,t;c_k)\Big)  \ep^*(\wh c_i)\ 
     =\ \dfrac d{dt} R_{ik}(t)
\end{eqnarray*}
$\text{in which we recall that}\quad R_{jk}(t)=\ds \sum_{\wh c_j=c_j^\pm}\theta(\wh c_j,t;c_k)\cdot \ep^*(\wh c_j).$\qed

{\noindent\em Proof of Theorem \ref{zeta}}. We have
  $$0=N_f(t) +Tr \FFF =N_f(t)  +Tr \RRR'\RRR^{-1}=N_f(t)  +\dfrac{\DDD'}\DDD. $$ 
  \qed

\subsection{Weighted lap function and proof of Theorem \ref{invertible}}

Let us consider the generating function of $g^n(\whx)$:
\REFEQN{generating} \begin{array}{ll} G(\whx,t) = \ds \sum_{n\ge 0} t^n g^n(\whx) \ \  \text{for } \whx\in \wh I \ \text{and  then}\\
 G(x,t) =\dfrac12\Big( G(x^\-,t)+G(x^\+,t) \Big)\ \  \text{when $a<x<b$}.\end{array} \ENDEQN

 Let $J=\eng{\whu,\whv}\subset ]a,b[ $ be an (open, closed or half-closed)
 interval or a point.  We define  the {\bf weighted lap function}\footnote{If $g_i\equiv 1$ the $G$-functions are $\dfrac1{1-t}$ and the function $L(J,t)$ is the 
generating function for the numbers of $(n+1)$-cylinders in $J$, and $L(]a,b[,t)$ has radius of convergence equal to $1/\rho_1$.}
\begin{eqnarray} L(J,t) &:=& \dfrac12\sum_{n\ge 0} t^{n}
  \sum_{\al\in Z_{n+1} }\sum_{\whx\in \wh \partial \al}  g^n(\whx)\chi_J(x). 
\end{eqnarray}
Repeating the calculation in our proof of the Main Kneading Identity
without the sign factors $s, \ep$ and $\sigma$, it follows easily that
\begin{eqnarray} L(J,t)  &=&
\sum_{j=1}^\ell\Big( 
   \sum_{p\ge 0,\, x\in \G_{p,c_j}}t^pg^p(x)\chi_J(x)\Big)
    \Big(\sum_{ \whc=c_j^\pm} \dfrac 12 \sum_{q\ge 0}t^q\cdot g^q(\whc)\Big)\\
    &=& \sum_{j=1}^\ell\g_{c_j,J}(t)\cdot G(c_j,t)  
    \label{gamma-identity}
\end{eqnarray}
In particular, for a one-point set $J=\{x\}$ we have simply
\REFEQN{one-point} L(\{x\},t)=\left\{\begin{array}{ll} t^p g^p(x) \cdot 
   G(c_i,t)  & \text{for } \ds x\in    \G_{p,c_i},\ p\ge 0, 
     \mystrut \ 1\le i\le \ell\\ 0 & \text{otherwise.}\end{array}\right.\ENDEQN

\REFLEM{analyticity} Fix any subinterval $J=\eng{\whu,\whv}$. The functions $G$, $\theta$, $\De_J \theta$, $\RRR_{jk}$ are all analytic functions of $t$ on the disc $\{ |t|<1/\rho_\infty\}$.
The kneading matrix is invertible when $|t|<1/\rho_1$.
The function $L(J,t)$ is meromorphic on $\{ |t|<1/\rho_\infty\}$ and analytic on $\{ |t|<1/\rho_1\}$.
\ENDLEM

\beginp 
 The first claim follows from the definition
of $\rho_\infty$ and the following estimates:

$$\forall\,\whx\in \whI,\quad 
   |G(\whx,t)|\le  \sum_{n\ge 0} |t|^n
     \|g^n\|_\infty < \infty \quad \text{for}\quad |t|<1/ \rho_\infty$$
     Similarly $\forall \, k$,  $$| \De_J \theta(\cdot,t;c_k)|\le  \sum_{n\ge 0} |t|^n
     \|g^n\|_\infty < \infty \quad \text{for}\quad |t|<1/ \rho_\infty.$$

To see that the kneading matrix is invertible when $|t|<1/\rho_1$
we use the relationship to the zeta function.
By Theorem \ref{zeta} we have $Z(t)\cdot \det \RRR(t)=1$,
where $$Z(t)=\exp \Big(\sum_{n\ge 1} \dfrac{N_n}n t^n\Big)$$
and each
$|N_n|\le 
\big\|g^n\big\|_1.$
So $Z(t)$ is analytic and non-zero for $|t|<1/\rho_1$ whence
$\RRR(t)$ is invertible for $|t|<1/\rho_1$. 

We have \REFEQN{norm-1}|L(J,t)|  \leq
    \sum_{n\geq 0} |t|^n \sum_{\al\in Z_{n+1}} \big|g^n_{|\al}\big| \leq 
    \sum_{n\geq 0} |t|^n \sum_{\al\in Z_{n}} \big|g^n_{|\al}\big|  (\ell+1)=(\ell+1)  \sum_{n\geq 0} |t|^n  \|g^n\|_1\ENDEQN
    which shows that $L(J,t)$ has radius of convergence at least $1/\rho_1$.

Using the MKI itself for the $\gamma$ factor in (\ref{gamma-identity})
we get:
    \begin{equation}
    L(J,t)=   \sum_{k=0}^\ell \De_J \theta(\cdot,t;c_k)
       \Big(\sum_{j=1}^\ell \RRR^{-1}(t)_{kj}\cdot G(c_j,t)\Big)
        \label{w-identity}
\end{equation}
The above identities are valid as formal power series but also
 when the functions involved are analytic and $\RRR(t)$ is invertible. 
As $1/\rho_1\le 1/ \rho_\infty$, so when $|t|<1/\rho_1 $, the identity \Ref{w-identity} is valid. 
\qed

{\noindent \em Proof of Theorem \ref{invertible}}
 
 The first two claims have already been proved in Lemma \ref{analyticity}.
 
 We proceed to prove the last claim.
When all $g_i$'s
are positive
and $t\geq 0$ we have
 \[ L(]a,b[,t)+ G(a^+,t) +G(b^-,t)=
    \sum_{n\geq 0} t^n \sum_{\al\in Z_{n+1}} g^n_{|\al} \geq 
    \sum_{n\geq 0} t^n  \|g^n\|_1 .\]
    By definition the RHS has radius of convergence equal to $1/\rho_1$.
Being a power-series with positive coefficients it follows that
the RHS diverges as $t\nearrow 1/{\rho_1}$.
    
Under the further assumption 
    $1/\rho_1 <1/ \rho_\infty$,  the functions $t\mapsto G(\whx,t)$, in particular
 $G(a^+,t)$ and $G(b^-,t)$, remain bounded at $t = 1/\rho_1$.
 So $L(]a,b[,t)$ must diverge as
$t\nearrow 1/{\rho_1}$.
Combining with (\ref{norm-1}) we know that the radius of convergence of  $L(]a,b[,t)$ is equal to $1/\rho_1$.
Now, the functions  $\Delta_J\theta$ and $G$  involved in \Ref{w-identity} remain bounded 
on $|t|\le 1/\rho_1$.
Letting $t\nearrow 1/{\rho_1}$ in \Ref{w-identity}  we conclude that $\RRR(t)$ must
be non-invertible at $t=1/{\rho_1}$.\qed

\section{Semi-conjugacies to piecewise linear models}\label{semi}
In this section we prove  Theorems \ref{conjugacy} and \ref{mainconjugacy}.

\REFLEM{forward}
Fix $J=\eng{\whu,\whv}\subset I_j= ]c_j^\+, c_{j+1}^\-[$.
We have for $k=0,\cdots,\ell$ and
$|t|<1/\rho_\infty$:
\begin{equation}\label{17}
\theta(\whv,t;c_k)-\theta(\whu,t;c_k) =
t \cdot s_jg_j \Big( \theta(f\whv,t;c_k)-\theta(f\whu,t;c_k)\Big)
\end{equation}
When also $|t| < 1/\rho_1$ we have for the weighted lap function~:
\begin{equation}
L(J,t)=  t\,g_j \cdot L(f_j J,t) \label{lap-weighted}
\end{equation}
\ENDLEM
\beginp Let us fix  $k\in \{0,\cdots, \ell\}$.
By definition, we have the following relation  for $\theta(\cdot,t;c_k)$
when applied to $\whx$ and $f\whx$:
$$\forall\,\whx\in \whI,\ \ \theta(\whx,t;c_k)=
   \sum_{m\ge 0}t^m [sg]^m(\whx)\cdot \sigma(f^m\whx, c_k)
= \sigma(\whx,c_k)+ t\cdot [sg](\whx)\cdot \theta(f\whx,t;c_k).$$ 
This implies \Ref{17} when restricting to $\whI_j$.
Now, 
$\De_J \theta(\cdot,t;c_k)=\theta(\whv,t;c_k)-\theta(\whu,t;c_k)$
and (as $f$ may reverse the orientation)
$\De_{f_jJ}\theta(\cdot,t;c_k),
= s_j\Big( \theta(f\whv,t;c_k)-\theta(f\whu,t;c_k)\Big)$, so
\begin{equation}
\De_J \theta(\cdot,t;c_k)=
 t g_j \De_{f_jJ}\theta(\cdot,t;c_k).
 \end{equation}

The result for $L(J,t)$ now follows by
linearity in equation (\ref{w-identity})
which is valid when $|t|<1/\rho_1$. \qed

{\noindent \em Proof of Theorem \ref{conjugacy}}.

We assume here that
all $g_i > 0$ and that $\rho_1>\rho_\infty$.
Fix $0<t< 1/\rho_1 < 1/\rho_\infty$.
Noting that $0< L(]a,b[,t)<+\infty$
we define our conjugating map
$\phi_t : \whI \rr \R$ by setting  \begin{equation}
  \phi_t(\whx) =  \frac{L(\eng{a^+,\whx},t)}{L(]a,b[,t)},
    \ \ \ \whx\in\whI.`
   \label{phi-def}
  \end{equation}
  Notice that $\phi_t$ maps point-germs to genuine real numbers.

 \begin{figure}
\epsfig{figure=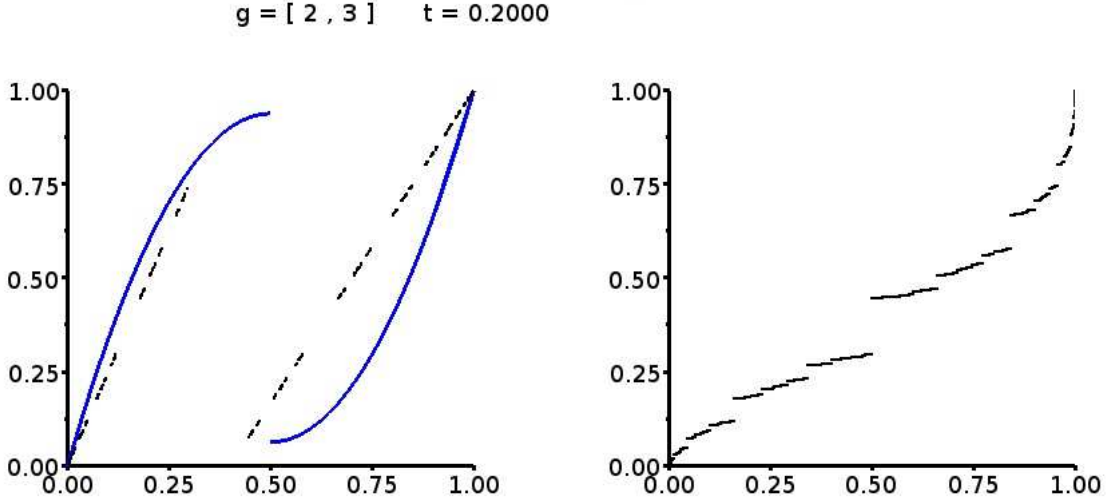,width=16cm}
\caption{Left: Example of a discontinuous map $f$ and the graph
restricted to $\Omega_t$ of its conjugated 
  map $\wt f_t$. 
  Right: The graph of $\phi_t$. 
  Here, $t=0.2 < 1/\rho_1=0.2684$.
  The ratio of slopes of the two branches is $3:2$
  (coming from the choice of weights).
  }
\label{discont figure}
 \end{figure}

{\noindent\bf Part A:}
  Using (\ref{lap-weighted})  we get for any 
$\whx_1,\whx_2\in \whI_j$ (the sign enters again) :
\begin{equation}
   \phi_t(\whx_2) - \phi_t(\whx_1) = 
       t s_j g_j \Big( \phi_t(f \whx_2) - \phi_t(f \whx_1) \Big).
         \label{phi-iterate}
\end{equation}
Similarly, we get by iterating this argument  for
$\whx_1,\whx_2\in \wh \al$ with $\al \in Z_n$ :
\begin{equation}
   \phi_t(\whx_2) - \phi_t(\whx_1) = 
       t^n \, s^n_{|\al} \; g^n_{|\al}
        \;\Big( \phi_t(f^n \whx_2) - \phi_t(f^n \whx_1) \Big).
         \label{alpha-iterate}
\end{equation}

When the $g_i$'s are non-negative, we clearly have 
$L(J,t)\geq 0$ for any interval $J$ so by set-additivity
with respect to $J$ it follows that
$\phi_t$ is monotone increasing and takes values in $[0,1]$.
Let $\Omega_t = \phi_t(\whI) \subset [0,1]$ and set
$\Omega_{t,i} = \phi_t(\whI_i)$. 
By monotonicity of $\phi_t$ the convex hull of $\Omega_{t,i}$
is precisely
 $\wt I_{t,i} = [ \phi_t(c_i^+), \phi_t(c_{i+1}^-)]$

Let now $a<x<b$. As all $g_i>0$,
by  \Ref{one-point}
\begin{equation}
\phi_t(x^+)-\phi_t(x^-)=\dfrac{L(\{x\},t)}{ L(]a,b[,t) }>0
    \label{cut pos}
    \end{equation}
precisely when $x$ is an interior cutting point or a pre-image
of such.  We have in particular $L(\{c_i\},t)\geq 1$ so that
$\sup \Omega_{t,i} < \inf \Omega_{t,i+1}$ and also
$\sup \wt I_{t,i} < \inf \wt I_{t,i+1}$, proving claim {\bf A.}

{\noindent\bf Part B:}
Given $y\in \Omega_t$ suppose that $y=\phi_t(\whx_1)=\phi_t(\whx_2)$
with $\whx_1<\whx_2$.  By the previous paragraph 
$\whx_1$ and $\whx_2$ must belong to the same  $\whI_j$.
(In fact they even belong to the same $n$-cylinder for all $n$).
So by the identity (\ref{phi-iterate}) we must have
$ \phi_t(f \whx_2) - \phi_t(f \whx_1) =0$.
This implies that there is
a well-defined map $\wtf_t : \Omega_t \rr \Omega_t$ given by:
\begin{equation}
   \wtf_t (y) :=  \phi_t (f \whx) , \ \ \ y=\phi_t(\whx) \in \Omega_t
\end{equation}
(since the value is independent of the choice 
of $\whx$ in the pre-image of $y$).

Equation (\ref{phi-iterate}) shows that the
conjugated map has 
(finite) slope $(ts_jg_j)^{-1}={s_j}/{tg_j}$ 
on each $\Omega_{t,j}=\phi(I_j)$ (if it is not reduced to a point).
The map $\wt f_{t,i}$ is defined to be this affine map extended
to $\wt I_{t,i}=[\phi_t(c_i^+),\phi_t(c_{i+1}^-)]$. \qed\\

{\noindent\bf Part C:}
The last part of the theorem is tricky due to
the fact that $\phi_t$ is neither continuous nor injective.
For an open interval $J=]u,v[ \subset ]a,b[$ we
will in the following use the short-hand notation:
  \begin{equation}
     \Xi_t(J) := [\phi_t (u^+), \phi_t(v^-) ]
  \end{equation}
For $0\leq i \leq \ell$ we define
$\wt I_t(i) = \wt I_{t,i} = \Xi_t(I_{i})$ 
and then recursively 
$
\wt I_t(i_0,\ldots, i_{n-1}) =  \wt I_{t,i_0} \cap \wt f_t^{-1} 
 \wt I_t(i_1,\ldots, i_{n-1})
 $
 which 
 is either empty, a point or
 a closed interval.
 We write $\wt Z_{t,n}$ for the collection of 
 non-empty sets of this form.
 They form a partition for the domain of definition of 
 $( \wt f_t )^n$. The maximal invariant domain for $\wt f_t$ is 
 the compact set
 $\ds \Omega_t= \bigcap_{n\geq 1} \left( \bigcup \wt Z_{t,n} \right) 
 \subset [0,1]$.
Our first goal is to exhibit a simple relationship between cylinders
and the above sets.

\REFLEM{phi-beta}  There is a bijection between
 $\alpha \in Z_n$ and $\wt \al \in \wt Z_{t,n}$
 given by $\wt \al = \Xi_t(\al)$.
\ENDLEM
Proof: 
For $n=1$ this is the very definition: $Z_1$ consists
of the intervals $\{I_i: 0\leq i\leq \ell\}$ and
  $\wt I_{t}(i)=\Xi_t(I_i) = [ \phi_t(c_i^+), \phi_t(c_{i+1}^-)]$.

When $J=]u,v[\subset I_i$ the definition of $\wt f_t$ shows that
$\wt f_t\; \Xi_t(J)= \wt f_{t,i} [\phi_t(u^+),\phi_t(v^-)]
 = [\phi_t(f_i u^+); \phi_t(f_i v^-)] = \Xi_t(f J)$.
For $\al\in Z_n$ this implies
$(\wt f_t)^k\; \Xi_t (\al) = \Xi_t(f^k \al)  \subset  \Xi_t(I_{i_k})$.
It follows by recursion that
$\Xi_t(\al) \subset \wt \al=\wt I_t(i_0,\ldots,i_n)$.

In order to show equality 
we proceed by induction in $n$. 
Let $\beta=(i_1 \ldots i_n)=]u_1,u_2[ \in  Z_n$  (with $a\leq u_1<u_2\leq b$).
Our induction hypothesis is that 
$\wt I_t(i_1,\ldots, i_n) = \Xi_t(\beta)= [\phi_t(u_1^+), \phi_t(u_2^-)]$.
Here, $u_1$ and $u_2$ are necessarily (eventual) cutting points
so by 
    (\ref{cut pos})
we have when $a<u_1$ and $u_2<b$, respectively :
\begin{equation}
\phi_t(u_1^-) < \phi_t(u_1^+)  
 \ \ \ \mbox{and} \  \ \ \ \phi_t(u_2^-) < \phi_t(u_2^+).
\label{strict ineq}
\end{equation}
Suppose that 
$\wt f_t\; \wt I_{t,i_0} = \Xi_t (f I_{i_0})$ intersects $\Xi_t(\beta)$
non-trivially
 and write $\wt \al = \wt I_t(i_0,\ldots,i_n)=[\xi_1,\xi_2]$.
We claim that
also $]v_1,v_2[  \equiv  f I_{i_0}$
intersects $\beta$. If this were not the case,
then e.g.\ $v_1<v_2\leq u_1<u_2$ in which case the first
strict inequality in 
(\ref{strict ineq}) 
shows that
$\phi_t(v_2^-)\leq \phi_t(u_1^-) < \phi_t(u_1^+)$
so that $\wt \al$ was empty in the
first place. 

Assume $s_{i_0}=+1$. We have
$\wt \al = [\xi_1,\xi_2]=
{\wt I}_{t,i_0} \cap \wt f_t^{-1} \Xi_t(\beta) $
(using the induction hypothesis)
and we write $\al = I_{i_0} \cap f^{-1} \beta = ]w_1,w_2[$
(which is non-empty as just shown).
We will calculate the left end points of  $\al$ and $\wt \al$.
We consider the two possibilities:
Either
$u_1 \leq v_1 (< u_2, v_2)$ or
$v_1 <  u_1 (< u_2, v_2)$.

In the first case $w_1^+=f_{i_0}^{-1} \max\{u_1^+, f c_{i_0}^+\} = c_{i_0}^+$ 
and since
$\phi_t(u_1^+) \leq \phi_t(v_1^+) = \wt f_t \phi_t(c_{i_0}^+)$
we get $\xi_1 = 
(\wt f_{t,i_0})^{-1} \max \{ \phi_t(u_1^+), \wt f_t  \phi_t(c_{i_0}^+) \} = 
\phi_t(c_{i_0}^+) = \phi_t(w_1^+)$.
In the second case, continuity and strict monotonicity of $f_i$ yields
a unique value $w_1 \in ]c_{i_0}, c_{i_0+1}[$ for which $fw_1=u_1$.
We have $\wt f_t \phi_t(w_1^+) = \phi_t(u_1^+) \leq \phi_t(f c_{i_0}^+) \leq 
  \wt f_t \phi_t(c_{i_0}^+)$ so again $\xi_1=\phi_t(w_1^+)$.

Thus in either case $\xi_1=\phi_t(w_1^+)$. Similarly, $\xi_2=\phi_t(w_2^-)$
thus implying $\wt \al = \Xi_t(\al)$ as we wanted to show.
If $s_i=-1$ some intervals change direction but the conclusion remains
the same.
\qed 

{\noindent \it  Returning now to the proof of Part C:}
Clearly $\phi_t(\wh I) \subset \Omega_t$. In order
to show surjectivity consider $\xi\in \Omega_t$.
Assume that $\wt f_t^k \in \wt I_{t,i_k}$, $k\geq 0$.
Then $\xi \in \wt \al_k = \wt I_t(i_0,\ldots,i_{k-1}) = \Xi_t(\al_k)$
for all $k$ (a nested sequence of intervals).
First, if $\xi$ is a boundary point of such an interval for some $k$ then
it is in the image of $\phi_t$ by the previous lemma. So assume
that $\xi$ 
is in the interior of $\wt \al_k = \Xi_t(\al_k)$ for all $k$.
Let $\al_k = ]u_k,v_k[$. 
Then $u_k\nearrow u_*$ and $v_k\searrow v_*$ with $u_* \leq v_*$.
None of the sequences are eventually constant.
Now, $\phi_t(u_k^+) \leq \xi \leq \phi_t(v_k^-)$ and
$0 \leq \phi_t(v_k^-)-\phi_t(u_k^+) \leq t^k g^k_{|\al_k}/L(]a,b[,t) \rr 0$
as $k\rr \infty$. For any $x\in [u^*,v^*]$  we conclude by
monotonicity of $\phi_t$ that
$\phi_t(x^+)=\phi_t(x^-) = \xi$. So $\phi_t : \wt I \rr \Omega_t$ 
is surjective.
\qed\\

In order to prove Theorem \ref{mainconjugacy}
we consider the limit
$t\nearrow 1/{\rho_1}$.
As the function $L(]a,b[,t)$ diverges the situation is a bit different.
By Lemma \ref{analyticity} the lap-function
  $L(]a,b[,t)$ is meromorphic
in  the disc  $\{|t|< 1/{\rho_\infty}\}$ and has a pole
of some order $m\ge 1$ at $t=1/{\rho_1}$. 
By positivity of $L(]a,b[,t)$ for $t>0$ there is $c>0$ so that
$$L(]a,b[,t)=\dfrac c{(1-{\rho_1}t)^m } + l.o.t.$$
For any interval $J\subset ]a,b[$ we  have
$0\le L(J,t)\le  L(]a,b[,t)$. An eventual  pole of $L(J,t)$ at 
$1/{\rho_1}$ is therefore of order at most $m$ 
so $\dfrac{L(J,t)}{L(]a,b[,t)}$ extends analytically
to $t=1/{\rho_1}$
(the singularity  is removable here).
We denote the limit
\REFEQN{measure}
\La (J):=\lim_{t\nearrow 1/{\rho_1}} 
\dfrac{L(J,t)}{L(]a,b[,t)} \in [0,1].\ENDEQN

 \begin{figure}
\epsfig{figure=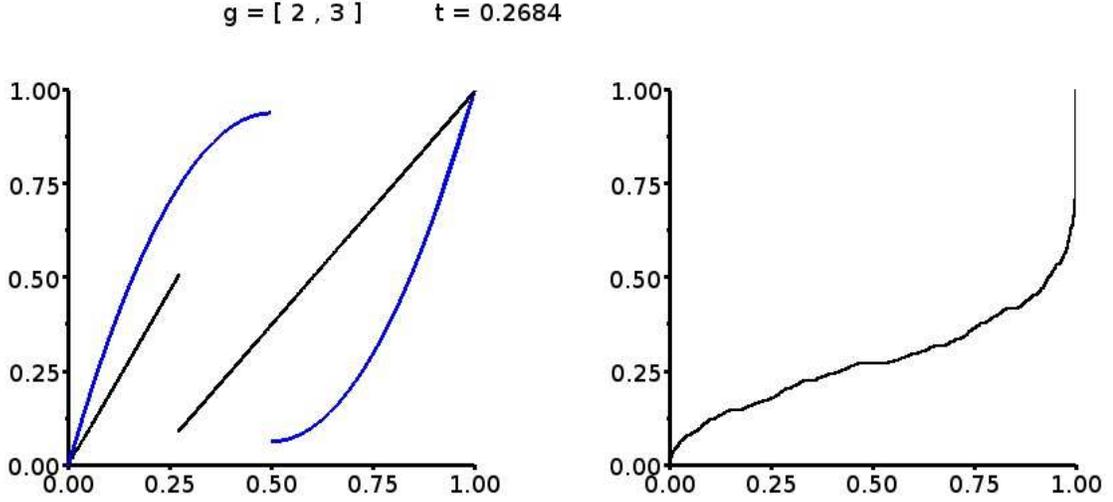,width=16cm}
\caption{The same example as before but at the critical value:
   $t=1/\rho_1=0.2684$. Note that $\Omega_t$ is
    no longer a Cantor set and that $\phi_t$ is continuous.
  }
\label{discont figure}
 \end{figure}

\REFLEM{n-cylinder}  We have the following properties for $\La$:
\begin{enumerate}
\item
For any $x\in I$ : $\La(\{x\})=0$. \mystrut
\item For all $\al\in Z_n$ : 
    $\ds \La(\al) = 
    \frac{1}{\rho_1^n}\;
    {g^n_{|\al}}\;
    \La(f^n\al) \;
    $. \mystrut
\item $ \delta_n = \ds \sup_{\al\in Z_n}\La(\al)\xrightarrow[n\to \infty]{}  0$.
   \mystrut
\end{enumerate}
\ENDLEM
\beginp 
The expression (\ref{one-point})
shows that the function $L(\{x\},t)$ is analytic on $\{|t|<1/\rho_\infty\}$ in particular  remains bounded on $\{|t|\le 1/\rho_1\}$. As $t\nearrow 1/\rho_1$,  the denominator $L(]a,b[,t)$ diverges, the first claim follows.

For $J\subset I_j$ for some $j$ we divide (\ref{lap-weighted})
by $L(]a,b[,t)$ and take the limit $t\nearrow 1/\rho_1$ to obtain:
\begin{equation}
    \La(J) = \frac{1}{\rho_1} g_j \La (f_j J) .
\end{equation}
In particular, for $\al=(i_0i_1\cdots i_{n-1})\in Z_n$ we have
$\al\in I_{i_0}$ so that 
$\La(\al)=\dfrac1{\rho_1} g_{i_0} \cdot\La(f\al).$
Iterating this we get the formula.

The last claim follows from
$$ \La(\al)=\dfrac1{\rho_1^n} \; {g^n_{|\al}}\;
    \La(f^n\al)\le \dfrac{\rho_\infty^n} {\rho_1^n}\  \xrightarrow[n\to \infty]{\rho_1>\rho_\infty}\  0.$$
 \qed

\REFLEM{phistarcont}
  The map $\phi : [a,b] \rr [0,1]$ defined by
   $$\phi(x)=\La(]a,x[), \ x\in [a,b]$$ is non-decreasing, continuous
   and surjective. One has for $x\in ]a,b[$:
   \begin{equation}
      \phi(x) 
          = \lim_{t\nearrow 1/\rho_1} \phi_t(x^-) 
          = \lim_{t\nearrow 1/\rho_1} \phi_t(x^+) 
      \label{phi-limit}
   \end{equation}
\ENDLEM
\beginp 
   Monotonicity follows from positivity and additivity of $\La(J)$, $J\subset ]a,b[$.
   Let $x\in [a,b]$ and $\ep>0$. Choose $n$ so that 
 $\ds \delta_n < \ep/2$ ($\delta_n$ from the previous lemma).  Either 
 $x$ is inside some $n$-cylinder or on the boundary
 of two such  cylinders. In any case, we may find
 at most two  n-cylinders $\alpha_1,\alpha_2$ with
 $\overline{\alpha_1} \cap \overline{\alpha_2} = \{x'\}$ so
 that $J=\alpha_1\cup\{x'\}\cup \alpha_2$ is an open neighborhood
 of $x$ and $\La(J) <\ep$. For $\delta>0$ small enough
 $\phi(x+\delta)-\phi(x-\delta) \leq \La(J)<\ep$.
 As $\phi(a)=\La(\emptyset)=0$ and $\phi(b)=1$ the map is surjective. 
 The first equality in (\ref{phi-limit})
 is essentially the definition of $\phi$ and
 the second follows from the continuity just shown.
 \qed
\\

We write $\wtc_i=\phi(c_i)$, $i=0,\cdots, \ell+1$ and let
$\wtS \subset S:=\{0,\cdots,\ell\}$ denote the
(possibly strict) subset of indices $i$ for which 
$0 < \wtc_{i+1} - \wtc_i = \La(]c_i,c_{i+1}[)$.
For $i\in \wtS$ we set
$\wtI_i=]\wtc_i,\wtc_{i+1}[$.

{\noindent \em 
Proof of Theorem \ref{mainconjugacy}.}

{\noindent\bf Part A:}
For $\whx_1,\whx_2\in I_j$,
taking the limit $t\nearrow 1/\rho_1$ in
the identity (\ref{phi-iterate})  
yields \begin{equation}
   \phi(\whx_2) - \phi(\whx_1) = 
       t s_j g_j \Big( \phi(f_j \whx_2) - \phi(f_j \whx_1) \Big).
\end{equation}
Continuity of $\phi$ and $f_j$ shows that this identity
is independent of the direction of the point-germs.
The affine map $\ds \wtf_j(y)= \whc_j + \frac{s_j}{t g_j} (y-\whc_j)$ then 
satisfies the required identity. 

{\noindent \bf 
Part B: }
Recall that $Z_n$ consists of the non-empty 
$n$-cylinder for $(I_i,f_i)_{i\in S}$.
 Let $\wtZ_n$ be the collection of non-empty
 open intervals of the form
$\wt\al=\Int\;  \phi(\al)$ where  $\al=(i_0 \cdots i_{n-1}) \in Z_n$. 
Here each $i_k\in \wtS$, $0\leq k< n$
(or else $\wt \al$ \`a fortiori empty) and 
$\wtf^k \wt \al \subset \wtI_{i_k}$. 
Therefore $\wt \al$ is contained in  an $n$-cylinder for the dynamical system
$(\wtI_i, \wtf_i)_{i\in \wtS}$.
We claim
that  $\wt \al$ is actually  equal to an  $n$-cylinder for that system and  $\wtZ_n$ is precisely the set of non-empty $n$-cylinders for
the same system. To see this note that
\begin{equation}\label{28}
1=\sum_{\al\in Z_n}\La(\al)
  =\sum_{\wt\al\in \wtZ_n}|\wt\al|,
  \end{equation}
where $|\cdot|$ denotes the length of intervals. 
There is no room for any other or any larger open cylinder.

Now, by Lemma \ref{n-cylinder}  we have
$|\wt\al|=\La(\al)=\dfrac{g^n_{|\al}}{\rho_1^n}
      \La(f^n\al) \le \dfrac{g^n_{|\al}}{\rho_1^n}$.
So using \Ref{28} we get
$$\rho_1^n=\sum_{\wt\al\in \wtZ_n} |\wt\al| \rho_1^n = \sum_{\wt\al\in \wtZ_n}\La(\al)\rho_1^n =  \sum_{\wt\al\in \wtZ_n} g^n_{|\al} 
  \La(f^n\al)\le \sum_{\wt\al\in \wtZ_n} g^n_{|\al}
   \le 
   \sum_{\al\in Z_n} g^n_{|\al}=\|g^n\|_1.
   $$
   So $\rho_1=\ds\limsup_{n\to\infty}\Big( \sum_{\wt\al\in \wtZ_n} g^n_{|\al}\Big)^{1/n}$.
The pressures of
 $(I_i,  f_i,g_i)_{i\in S}$ 
and $(\wtI_i,  \wtf_i,g_i)_{i\in \wtS}$ are 
therefore the same.

{\noindent \bf Part C: }
We assume here that $f$ extends to a continuous map of $[a,b]$.
When $J\subset [a,b]$ is an interval then 
$fJ \setminus \bigcup_i f(J\cap I_i)$ consists of a finite
number of points.  By Lemma \ref{n-cylinder} this
set difference has zero mass. By the same lemma
we get:
$\La(J)\ds=\sum_{i=0}^{\ell}  \dfrac1{\rho_1} g_i  \La(f_i(J\cap I_i))$.
Thus,
\[\Big(\min_i g_i\dfrac1{\rho_1}\Big)\La(fJ)\le \La(J)\le
     \Big(\sum_i g_i\dfrac1{\rho_1}\Big)\La(fJ).\]
In particular $\La(J)=0\Longleftrightarrow \La(fJ)=0$.
(Note, however, that $\La(J)>0$ does not imply $\La(f^{-1}J)>0$ as the 
latter set might be empty). 

Let us write $x\sim x'$ if $\phi(x)=\phi(x')$,

When $x,x'\in I$ and $x\sim x'$ then $\La([x,x'])=0$ 
so also $\La(f[x,x'])=0$. 
As we have assumed $f$ continuous,
$f([x,x'])$ is connected and contains $f(x)$, $f(x')$. 
Therefore,
$\phi(f(x))=\phi(f(x'))$, i.e. $f(x)\sim f(x')$. 
For $y\in [0,1]$, we may thus define
$\wtf(y)=\phi(f(x))$ with $x\in \phi^{-1}(y)$
(independent of the choice of $x$).
Then for every $x \in [a,b]$ we have:
\[\wtf( \phi(x))=\phi(f(x))\]
The same argument also shows that for any two $x,x'\in I$
we have 
$$\ds \Big|\wtf(\phi(x))-\wtf(\phi(x'))\Big| \leq \max_{i\in \wtS} \frac{\rho_1}{g_i}
   |\phi(x)-\phi(x')|$$ so $\wtf$ is a continuous endomorphism of $[0,1]$.
\qed

Remark: The set $\wtS$ may depend upon the weights $g_i$.
If, however, $f$ is transitive then $\wtS=S$ for any choice
of non-zero weights and $\wtZ_n=Z_n$ for all $n$.
We leave the exercise to the reader.

\appendix

\section{Geometry of the weight function $\omega(\al)$}\label{count}

Fix $n\ge 1$ and an $n$-cylinder $\al\in Z_n$. Recall that we have associated a weight 
$$\omega(\al)=-\; g^n_{|\al}  \sum_{\whx\in \wh \partial \al } \sigma(f^{n}\whx,x)\cdot \ep(f^n\whx).$$
Set 
$$\pi(\al):=-\sum_{\whx\in \wh \partial \al } \sigma(f^{n}\whx,x)\cdot \ep(f^n\whx).$$
This quantity depends only on the boundary values and
their  positions relative to the diagonal. Let  $h$ be an affine map on $\al$ coinciding with $f^n$ on the boundary.

\REFLEM{weight} $ \pi(\al)\ds=-\sum_{\whx\in \wh \partial \al } \sigma(h(\whx),x)\cdot \ep(h(\whx))$. And,\\
$\bullet$ $\pi(\al)=-1$ if $0<\text{\rm slope}(h)\le 1$ and $\overline{h(\al)}$ touches the diagonal\\
$\bullet$ $\pi(\al)=1$ if $h(\al)$ transverses the diagonal with slope either $>1$ or $<0$.\\
$\bullet$ $\pi(\al)=0$ in all other cases, namely \\
\mbox{} \quad either $\overline{h(\al)}$ does not touch the diagonal\\
\mbox{} \quad or $\overline{h(\al)}$ touches the diagonal at one end only, with slope $>1$ or $<0$.
\ENDLEM
\begin{figure}
\scalebox{1} 
{
\begin{pspicture}(0,-2.54)(12.301894,2.54)
\psdots[dotsize=0.12](0.06,-2.46)
\psline[linewidth=0.02cm](0.06,-2.46)(4.86,2.34)
\psline[linewidth=0.04cm,arrowsize=0.05291667cm 2.0,arrowlength=1.4,arrowinset=0.4]{->}(0.04,-2.46)(0.06,2.5)
\psline[linewidth=0.04cm](0.98,0.32)(1.34,-0.16)
\psline[linewidth=0.04cm](1.58,-0.18)(2.02,0.24)
\psline[linewidth=0.04cm](2.14,-1.0)(2.42,-0.4)
\psline[linewidth=0.04cm](2.58,-0.34)(3.08,-0.96)
\psline[linewidth=0.04cm](4.12,1.6)(4.36,2.4)
\psline[linewidth=0.04cm](3.66,1.92)(4.0,1.44)
\psline[linewidth=0.04cm](3.48,0.96)(3.8,0.2)
\psline[linewidth=0.04cm,arrowsize=0.05291667cm 2.0,arrowlength=1.4,arrowinset=0.4]{->}(0.06,-2.48)(5.62,-2.46)
\psline[linewidth=0.04cm,arrowsize=0.05291667cm 2.0,arrowlength=1.4,arrowinset=0.4]{->}(6.68,-2.46)(12.24,-2.44)
\psline[linewidth=0.04cm,arrowsize=0.05291667cm 2.0,arrowlength=1.4,arrowinset=0.4]{->}(6.64,-2.44)(6.66,2.52)
\psdots[dotsize=0.12](6.64,-2.46)
\psline[linewidth=0.02cm](6.66,-2.44)(11.46,2.36)
\psline[linewidth=0.04cm](7.86,-0.72)(8.62,-0.46)
\psline[linewidth=0.04cm](7.5,-1.6)(8.12,-1.46)
\psline[linewidth=0.04cm](8.64,-0.14)(9.72,0.32)
\psline[linewidth=0.04cm](10.98,2.4)(11.54,1.78)
\psline[linewidth=0.04cm](10.34,0.84)(10.86,2.08)
\psline[linewidth=0.06cm](6.78,-2.28)(7.34,-1.72)
\psline[linewidth=0.04cm,linestyle=dashed,dash=0.16cm 0.16cm](9.98,2.3)(9.98,-2.32)
\usefont{T1}{ptm}{m}{n}
\rput(11.351455,0.525){$\pi(\al)=+1$}
\usefont{T1}{ptm}{m}{n}
\rput(7.781455,0.345){$\pi(\al)=-1$}
\usefont{T1}{ptm}{m}{n}
\rput(3.701455,-1.935){$\pi(\al)=0$}
\psline[linewidth=0.04cm](3.1,0.24)(3.38,0.84)
\end{pspicture} 
}\caption{fixed points counting}\label{Fix}\end{figure}
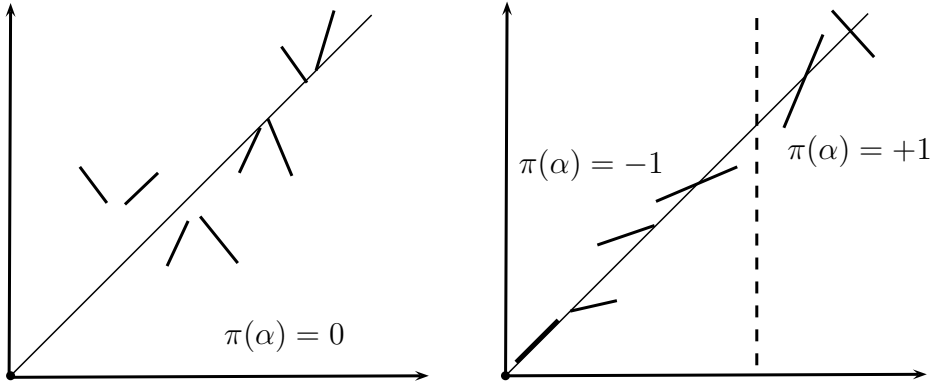

\beginp Since $f^n|_\al$ is a continuous strictly monotone map, we have  $\sigma(f^n\whx,x)=\sigma(h(\whx),x)$ and $\ep(f^n\whx)=\ep(h(\whx) )$ at 
the two ends of $\al$. So we can replace $f^n$ by $h$ in $\pi(\al)$.

Extend $h$ continuously to the boundary points. For $\whx$ a boundary germ.
We check case by case  the value of 
$-\sigma(h(\whx),x)\cdot\ep(h(\whx) )$ 

is  $\dfrac12$ if $h(x) <x$ and $h(\whx)>h(x)$, or if $h(x) >x$ and $h(\whx)<h(x)$

is $-\dfrac12$ if $h(x)=x$, or $h(x) >x$ and $h(\whx)>h(x)$, or if $h(x) <x$ and $h(\whx)<h(x)$.

Adding the values at the two ends, we get the lemma. \qed

Notice that if $f$ is expanding then $\pi(\al)\ge 0$ for all $n$ and all $\al\in Z_n$. So in that case $N_n\ge 0$ for all $n$.

\section{Relation between $\det\RRR$, $\det \BBB$ and Milnor-Thurston's kneading determinant}
We relate here our definition of the kneading determinant to that of Milnor-Thurston, modified by adding weights. 
Set $\theta(\wh x,t;c_{\ell+1}):=-\theta(\wh x,t;c_0)$

\REFLEM{key}(Minor-Thurston) We have  $\ds
\sum_{k=0}^{\ell} (1-t s_kg_k)(\theta(\whx,t;c_k)-\theta(\whx,t;c_{k+1}))\equiv 1.$\ENDLEM

For $k=0,\cdots, \ell$, set $I_k=]c_k, c_{k+1}[$.  Note that $\sigma(\whx, c_k)-\sigma(\whx, c_{k+1})= \chi_{I_k}(\wh x)$. 
Set $$\eta(\whx,t;I_k):=\theta(\whx,t;c_k)-\theta(\whx,t;c_{k+1})=\sum_{m\ge 0} t^m [sg]^m(\whx)\chi_{I_k}(\whx).$$
$$\De_{c_i}\eta(\cdot,t;I_k):=\eta(c_i^\+,t;,I_k)-\eta(c_i^\-,t;I_k), \ i=1,\cdots, \ell$$
And define the {\bf Milnor-Thurston kneading matrix} $\ell \times ( \ell+1)$ matrix $$\NNN(t)=\Big(\De_{c_i}\eta(\cdot,t;I_k)\Big)_{ i=1,\cdots,\ell,\ k=0,1,\cdots, \ell }\ .$$
Denote by $D_j$ the determinant of $\NNN(t)$ after deleting the $j$-th column.

\REFLEM{basic}(Milnor-Thurston) The quantity
$\dfrac{(-1)^jD_j}{1-s_jg_jt} =: D_{MT}(t)$ is independent of $j$ and is called the {\bf Milnor-Thurston kneading determinant}.\ENDLEM

\beginp  
Let ${\bf v}= \lv 1-s_0g_0t\\  \vdots \\  1-s_\ell g_\ell t \rv$. By Lemma \ref{key}, 
$\Big(\eta(\whx, t;I_0),\cdots, \eta(\wh x,t; I_\ell)\Big){\bf v}=1.$ So $ \bv$ is a kernel vector of $\NNN(t)$.
Define an augmented kneading matrix $\AAA(t)$ by adding a line vector $(\dfrac1{1-s_0g_0t}, \cdots , \dfrac1{1-s_\ell g_\ell t})$ on top of $\NNN(t)$. 

Then $\AAA\bv=\lv \ell+1\\ 0\\ \vdots \\ 0\rv.$ By Cramer's solution form  $1-s_jg_j t=(\ell+1) \dfrac{(-1)^ijD_j}{\det \AAA}$ and therefore
 $\dfrac{\det\AAA}{\ell+1} =\dfrac{(-1)^jD_j}{1-s_jg_jt} $  for  all $j=0,\cdots,\ell$.\qed

\REFLEM{relation} Setting  $H(t):=1-t\dfrac{s_0g_0+s_\ell g_\ell}2$  we have  $$ D_{MT}(t) =\det(\RRR(t))\quad \text{and}\quad H(t)\cdot \det \RRR(t)=\det \BBB(t).$$ \ENDLEM
\beginp
Set $\kappa_i(t)=\dfrac t2(s_{i-1}g_{i-1}-s_ig_i)$, $i=1,\cdots, \ell$. 
Grouping the terms about $c_0$ and $c_{\ell+1}$ in Lemma \ref{key} we get   
$$2H(t)\cdot \theta(\whx,t;c_0)+\sum_{i=1}^{\ell}2\kappa_i\cdot\theta(\whx,t;c_i)\equiv 1$$
It follows that
$$\lv H(t) & \kappa_1(t) & \kappa_2(t) &\cdots & \kappa_\ell(t) \\
0 & \\
0\\
\vdots & & Id \\
0\rv\cdot  \RRR(t)=\lv 1 & 0 & 0 &\cdots & 0\\ 1 \\ 1\\
\vdots&& \BBB \\
1\rv$$
Therefore $H(t)\cdot \det \RRR(t)=\det \BBB(t)$. For the matrix $\AAA$ defined above, 
$$ \AAA(t)\cdot \dfrac12 \lv  1-s_0g_0t & -1 & -1 &\cdots &-1\\
1-s_1g_1 t & 1 & -1 & \cdots & -1 \\
1-s_2 g_2 t & 1 & 1 & \cdots & -1 \\
\vdots\\
1-s_\ell g_\ell t & 1 & 1 &\cdots & 1 \rv
=\lv\frac{\ell+1}2 & * & * &\cdots & *\\ 0 \\ 0\\
\vdots&& \BBB \\
0\rv$$
 In the second matrix  on the left hand side add the last line   to every other line one gets \\
$\dfrac{\ell+1}2D_{MT}(t) \cdot H(t)=\dfrac{\det \AAA}2 \cdot  H(t)=\dfrac{\ell+1}2 \det\BBB.$

Combining with Lemma \ref{relation} we get  $D_{MT}(t)\cdot H(t) =\det\BBB= \det\RRR \cdot H(t)$. Therefore $\det \RRR(t)=D_{MT}(t)$.
 \qed

\REFCOR{weight-one} If all the weights $g_i$ are equal to $1$, all three determinants $D_{MT}$, $\det\RRR$, $\det\BBB$ have the same zeros
in $\D$.\ENDCOR

\beginp In this case   $H(t)=1-\dfrac t2(s_0+s_\ell)=1$ or $1-t$ so $H(t)$ has no zeros in $\{|t|<1/\rho_\infty\}=\D$.\qed

\section{The first zero of $\det \BBB$ may not  correspond to the pressure}\label{example}

We have shown in  Theorem \ref{invertible}, Point 3,  that the first zero of $\det \RRR$ corresponds to the pressure. And in case all the weights $g_i$ are $1$, one can also use the
first zero of $\det \BBB$ (Corollary \ref{weight-one}). 
This need not, however, be  true with more general weights.
Here is a counter example.

Let $I=[a,b]=[0,3]$, $I_0=]0,1[$, $I_1=]1,2[$, $I_2=]2,3[$.
$$f(x)=\left\{\begin{array}{ll} 2x & 0\le x \le 1\\
2-2(x-1) & 1\le x\le 2 \\ 2(x-2) & 2\le x \le 3 \end{array}\right.$$

Let us assign weights $g_0=g_1=1$ and $g_2=M$.

Note that $f(I_2)=[0,2]$ and $f: [0,2]\to [0,2]$ is the full tent map. There is no periodic points in $I_2$. 
Using Lemma \ref{weight} and the definition one obtains
$$Z(t)=\exp\Big(\ds \sum_{n\ge 1} \dfrac {t^n}{n} 2^n\Big) =(1-2t)^{-1}.$$
So by Lemma \ref{relation} and Theorem \ref{zeta} we have $D_{MT}(t)=\det\RRR(t) =\dfrac1{Z(t)}=1-2t$. The first zero being $1/2$ one obtains  that the pressure is $\log 2$ (this pressure can also be computed directly). It is easily seen that the topological entropy
is also  $\log 2$. 

On the other hand, $H(t) =1-\dfrac t2(s_0g_0+s_2g_2)=1-\dfrac t2(1+M)$. So by Lemma \ref{relation} again
$$\det\BBB(t)=H(t)\det\RRR(t) =\Big(1-\dfrac t2(1+M)\Big)(1-2t).$$
If $M>3$, then $\det\BBB(t)$ has a 'spurious' zero at $\dfrac2{1+M}$ smaller than $\dfrac 12$.  

So  the first positive zero of $\BBB(t)$ does not correspond to the pressure in this case. By increasing $M$, one can make this first zero arbitrarily small
without changing the pressure.

Hans Henrik RUGH,      
Bâtiment 425, Faculté des Sciences d'Orsay, Université Paris-Sud, 91405 Orsay Cedex, France

TAN Lei, 
Faculté des sciences,
LAREMA, Université d'Angers,
2 Boulevard Lavoisier, 49045 Angers cedex,
France 
\end{document}